\documentclass[a4paper,10pt]{article}
\usepackage[utf8x]{inputenc}
\usepackage{amssymb}

\title{On the existence of Lagrangian shadows of ample algebraic divisors. }
\author{Nikolay A. Tyurin\footnote{BLTPh JINR (Dubna) and  NRU HSE (Moscow), {\bf ntyurin@theor.jinr.ru}. {\it The author  was supported by RSF grant, project 14-21-00053 dated 11.08.14.}}}

\begin{document}

\maketitle

\begin{abstract} In a previous paper we introduce the notion of lagrangian shadows for ample divisors in algebraic varieties.
In this paper we present a condition of their existence.
\end{abstract}

In preprint [1] we presented the main constructions of Special Bohr - Sommerfeld geometry, a new programme which leads to the appearence of  certain  moduli spaces
in the framework of symplectic geometry. In preprint [2] this construction was specialized to the case of simply connected compact algebraic varieties which can be regarded
as symplectic manifolds equipped with compatible integrable complex structures. For any very ample divisor $D \subset X$ of an algebraic variety $X$ one can find a lagrangian shadow
$Sh^{Lag}(D) \subset X$ given by a finite number of (possible singular) lagrangian submanifolds in $X$ where the corresponding to $D$ Kahler form plays the role
of the symplectic form. At the same time this lagrangian shadow can be empty, see examples in [2].

Since all the arguments in the construction of lagrangian shadows follow the Andreotti - Frankel approach in the proof of the hard Lefschetz theorem [3], we can extend the material
of [2] adding the following proposition:

{\bf Theorem.} {\it The lagrangian shadow of the orientable type  of an ample divisor $D$ in a simply connected compact algebraic variety $X$ of complex dimension $n$ is non empty if and only if
the group $H_n(X \backslash D, \mathbb{Z})$ is non trivial.}

Before claiming this one should modify the notion of lagrangian shadows discussing the types of singularities which one allows for $Sh^{Lag}(D)$. By the very construction, see [2],
this set is formed by finite trajectories of the gradient flow of  smooth function $\phi_D = - \rm{ln} \vert h_D \vert$, where $h_D$ is the holomorphic section defined up to constant
by the condition $(h_D)_0 = D$, on $X \backslash D$. In any case $Sh^{Lag}$ carries the structure of CW - complex, therefore if we say that $Sh^{Lag}$ is a CW -complex,
the Theorem above turns to be a tautology.  Now instead of counting irreducible components of $Sh^{Lag}$ as it was proposed in [2] or Special Bohr - Sommerfeld cycles
as in [1], we just say that the rank of $H_n(X \backslash D, \mathbb{Z})$ is the number of the components and understand any generator of $H_n(X \backslash D, \mathbb{Z})$ as a single
(and possible singular, of course) component of the shadow. 

Below we prove Theorem and discuss some other properites of Lagrangian shadows.

Let $(M, \omega)$ be a compact simply connected symplectic manifold of real dimension $2n$ with integer symplectic form $\omega$, so $\int_{\Sigma} \omega \in \mathbb{Z}$ for any closed orientable 2 - dimensional
submanifold $\Sigma \subset M$. A submanifold $S \subset M$ is called lagrangian iff it is $n$ - dimensional and $\omega|_S $ vanishes. In this paper we
consider the case of compact {\it orientable} lagrangian submanifolds. Let $(L, a)$ be the prequantization data:
 line bundle with $c_1(L) = [\omega]$ and hermitian connection s.t. $F_a = 2 \pi i \omega$. A lagrangian submanifold
$S \subset M$ is Bohr - Sommerfeld of level $k$ (or $BS_k$ for short) iff the restriction $(L^k, a_k)|_S$ admits a covariantly constant section. A $BS_k$ lagrangian submanifold
$S \subset M$ is special w.r.t. a section $h \in \Gamma (M, L^k)$ iff  $h$ doesn't vanish on $S$ and 
$$
h|_S = e^{ic} f \sigma_S
$$
where $c$ is a real constant, $f$ is a real strictly positive function and $\sigma_S$ is a covariantly constant section of $(L^k, a_k)|_S$. The details can be found in [1].
As well there one shows how the definitions can be reformulated to work not only for smooth lagrangian submanifold but for lagrangian cycles, given by lagrangian submanifolds with 
prescribed types of singularities.

Suppose that $(M, \omega)$ admits a compatible complex structure $I$ which is integrable. Then $(M, \omega, I)$ can be regarded as an algebraic variety with Kahler metric of the Hodge type,
see [GH]. It is well known that in this case the space $\Gamma (M, L^k)$ contains a finite  dimensional subspace of holomorphic sections $H^0(M, L^k)$; for every section
$h \in H^0(M, L^k)$ the zeroset $(h)_0 = D \subset M$ is a complex $n-1$ - dimensional submanifold, which we  call divisor following algebro - geometrical tradition.

As it was pointed out in [1] the strategy of the search of Bohr - Sommerfeld lagrangian cycle which is special w.r.t. a holomorphic section $h_D$, corresponding to divisor $D$, is the following.
One takes the plurisubharmonic function $\phi_D = - \rm{ln} \vert h_D \vert$ correctly defined on the complement $M \backslash D$, then one takes the finite critical points
and finite trajectories of the gradient flow, which join these finite critical points,  and form the union of all these trajectories denoting this as $B_D \subset
M \backslash D$. Every $n$ - dimensional subset of this set  which is a compact possible singular submanifold   must be special Bohr - Sommerfeld lagrangian, see [1].

By the very construction $B_D$ is formed by isotropical cells, it carries the structure of CW - complex and its topological type is exactly the same as $M \backslash D$.
Indeed, $M \backslash D$ is $B_D$ plus unbounded trajectories of the gradient flow of $\phi_D$ therefore every unbounded trajectory can be contracted to the begining point
which is a finite critical point belonging to $B_D$. Thus $B_D$ is the retract of $M \backslash D$, which implies the statement. 

Now  since $B_D$ can contains open $n$ - dimensional component we exclude all these bad components getting a union of compact $n$ - dimensional cycles
which are our special Bohr - Sommerfeld cycles we are interested in. 

In [2] this possible reducible possible singular isotropical cycle $B_D$ was called lagrangian shadow of ample divisor $D \subset M$ for the case when
$B_D$ itself is a union of compact  lagrangian submanifold. We denote these shadows as $Sh^{Lag}(D)$.  Moreover in the examples presented in [2] these shadows were smooth.

{\bf Remark.} In the examples from [1], [2] one can see, that for small levels the Lagrangian shadows are smooth; at the same time if we go up increasing the level $k$ then
the Lagrangian shadows turn to be singular. As an example we again take $M = \mathbb{C} \mathbb{P}^1$: for the level $k=2$ we always have smooth Lagrangian shadow, but for
the level $k = 3$ in general the Lagrangian shadows are not smooth being isomorphic to figure $\Theta$. Topologically it is the buquet of two circles, and it is possible to
find a section of ${\cal O}(3)$ which corresponds to pair of disjoint loops, but this case is very far to be generic. The same happens for the case of complex 2- dimensional quadric:
for ${\cal O}(1,1)$ lagrangian shadows in general are smooth and it is still true for ${\cal O}(2,2)$, but not for higher levels. It seems that certain thresholds exist and these ones
are somehow related to the degree of the anticanonical bundle in the case of Fano varieties.

  The Theorem above can help in the solution of different problems. First of all it indicates emptyness of Lagrangian shadows for certain cases. For example in
[2] we have claimed that in 2 - dimensional complex quadric $M = Q \subset \mathbb{C} \mathbb{P}^3$ the ample divisor $D \subset Q$ given by the plane section 
admits non empty Lagrangian shadow iff it is irreducible (so the plane is not tangent to $Q$). Now applying the Theorem above we can see it clearly: if $D$ is reducible then
it is the union of two lines therefore $Q \backslash D$ for this case is contractible being homotopic to the product $\mathbb{C} \times \mathbb{C}$. Another example, where the Theorem
works ``positively'': let $M$ be the projective plane $\mathbb{C} \mathbb{P}^2$ and $L = {\cal O}(3)$, then the Lagrangian shadow of a smooth elliptic curve $D \in \vert 3l \vert$
is non empty and must be homotopic to a 2 - torus. The same arguments show that in the previous example with $Q \subset \mathbb{C} \mathbb{P}^3$ an irreducible divisor
from the anticanonical linear system has non trivial Lagrangian shadow, homotopic to 2 - torus.

Another possible application of the Theorem is in the problem of much more general type. It was conjectured  that a smooth lagrangian torus which is Bohr - Sommerfeld of level one or two
in the projective plane $\mathbb{C} \mathbb{P}^2$ with the standard  Kahler form does not exist (recall that in the standard Clifford fibration the minimal Bohr - Sommerfeld level
is 3 for smooth lagrangian 2 - tori). Suppose that we can prove that if $S \subset \mathbb{C} \mathbb{P}^2$ is a Bohr - Sommerfeld lagrangian torus of level $k$ then it can be
deformed by a Hamiltonian isotopy to special with respect to a holomorphic section lagrangian torus. Then we could apply the Theorem above and get the contradiction since the group
 $H_2(\mathbb{C} \mathbb{P}^2 \backslash D, \mathbb{Z})$ is trivial when $D$ is a projective line or  a conic.

At the end of this small note we would like to present a very natural question about lagrangian shadows of ample algebraic divisors. Namely suppose that an ample
algebraic divisor $D \subset X$ in a simply connected smooth compact algebraic variety $X$ admits certain smooth connected  lagrangian shadow $S_D = Sh^{Lag}(D) 
\subset X$. Then one can ask about the Floer cohomology of $S_D$ whether or not $FH^i(S_D, S_D; \mathbb{Z})$ can be expressed in terms of $X$ and $D$.
It seems that it should be the case since $S_D$ is completely defined by $X$ and $D$. 

The work on this subject is in progress.

 $$$$
{\bf References}

[1] Nik. A. Tyurin, ``Special Bohr - Sommerfeld geometry'',  arXiv:1508.06804;

[2] Nik. A. Tyurin, ``Lagrangian shadows of ample algebraic divisors'', arXiv: 1601.05974;

[3] A. Andreotti, Th. Frankel, ``The Lefschetz theorem on hyperplane sections'', Ann. of Math.,  Vol. 69, No. 3 (May, 1959), pp. 713-717.

\end{document}